\documentclass[preprint,11pt]{article}

\usepackage{graphicx,epsf,amsbsy,amsmath,amssymb}
\usepackage[latin1]{inputenc}

\font\bbfnt=msbm10
\def\bbR{\mbox{\bbfnt R}}
\newcommand{\mb}[1]{\mbox{\bfseries \itshape #1}}

\begin{document}

\title{Flow curvature manifolds \\
for shaping chaotic attractors: \\
{\sc i} Rössler-like systems}

\author{Jean-Marc Ginoux,\\ Laboratoire {\sc Protee}, I.U.T. de Toulon,\\ Université du Sud, BP 20132
F-83957 La Garde Cedex, France,\\Christophe Letellier,\\ CORIA UMR 6614, Universit\'e de Rouen,\\ BP 12 F-76801 Saint-Etienne du Rouvray cedex, France}


\maketitle

\begin{abstract}
Poincar\'e recognized that phase portraits are mainly structured around fixed points. Nevertheless,
the knowledge of fixed points and their properties is not sufficient to determine the whole
structure of chaotic attractors. In order to understand how chaotic attractors are shaped by singular
sets of the differential equations governing the dynamics,
flow curvature manifolds are computed. We show that the time dependent components of such manifolds
structure R\"ossler-like chaotic attractors and may explain some limitation in the development
of chaotic regimes.
\end{abstract}

\date{{\bf Keywords}: Chaos topology ; Flow curvature manifold.}

\section{Introduction}

Since the recognition of the importance of chaotic attractors in the
description of physical phenomena \cite{Rue71,Gol75,Ros76a,Hak75},
interest in developing techniques to characterize chaotic behaviours has lead
to many different approaches that can be roughly classified into i) a
statistical approach related to the ergodic theory \cite{Eck85,Aba93b} and
ii) a topological approach \cite{Gil02}. The characterization of chaotic
behaviours is a rather mature problem, at least for the three-dimensional
cases. In particular, the different types of chaos that can be encountered in
three-dimensional phase spaces are now well documented \cite{Gil02,Tsa03,Let06a}.
In spite of that, little has been said about the algebraic structure
that the differential equations must have for producing chaos. It is known since
Poincaré's early works that equations
describing chaotic flows must be nonlinear, non-integrable and at
least three-dimensional, according to the Poincaré-Bendixson theorem
\cite{Poi90,Ben01}.

These conditions are necessary but not sufficient
to produce chaos. Recently, it has been proved that quadratic
systems of ordinary differential equations, with a total of four terms
on the right-hand side, cannot produce chaotic attractors
\cite{Fu97}. In other words, a fifth term is required
to produce a chaotic attractor. From this point of
view, the minimal algebraic structure of a set of three ordinary
differential equations that produce a chaotic attractor corresponds to
four linear terms and one nonlinear term in the right-hand side (see
\cite{Spr00} for a review of investigations to discover simpler examples of
chaotic flows than the Lorenz and Rössler systems). Sprott was able to identify two
minimal equivalent chaotic flows \cite{Spr97}, whereas Malasoma \cite{Mal02} found
seven new examples of such minimal flows. These nine chaotic systems can be grouped
into two distinct classes
\cite{Mal02}. Nevertheless, nothing is said about the topology of their chaotic solutions.

Indeed, although fixed points have a prominent role in structuring the phase portrait, the whole shape of the attractor cannot be deduced
from them. Recently, it has been established \cite{Gin06,Gin08} that local metric
properties of chaotic attractors like the {\it curvature of the flow} can be analytically computed.
The set of points where the curvature vanishes defines the so-called {\it flow curvature manifold}
for which the invariance under the flow was proved by the Darboux theorem \cite{Gin06,Gin08,Dar78}. The aim
of this paper is to show that the time dependent component of the flow curvature manifold plays an
important role in the structure of chaotic attractors. The subsequent part of this paper is
organized as follows. In section \ref{flocur}, the procedure to compute flow curvature manifold
is detailed and its topology in the neighborhood of the fixed points is described. Section
\ref{Rosso} is devoted to explicit examples of many R\"ossler-like attractors. Section \ref{conc}
gives a conclusion.

\newpage

\section{Flow curvature manifold for 3D linear flows}
\label{flocur}

Let us consider the set of differential equations
\begin{equation}
  \label{sysdyn}
  \dot{\mb{X}} = \frac{{\rm d} \mb{X}}{{\rm d}t} = \mb{F} (\mb{X})
\end{equation}
where $\dot{\mb{X}}$ is the velocity vector. The state vector is such that
\begin{equation}
  \mb{X} = \left[ \displaystyle x_1, x_2, ..., x_n \right]^t \in E \subset \bbR^n
\end{equation}
and
\begin{equation}
  \mb{F} (\mb{X}) =
  \left[ \displaystyle F_1 (\mb{X}), F_2 (\mb{X}), ..., F_n (\mb{X}) \right]^t \in E \subset \bbR^n \, .
\end{equation}
The vector field $\mb{F} (\mb{X})$ is defined in a subspace $E$ in which its components $F_i$ are supposed to be continuous
and infinitely differentiable with respect to all $x_i$ and $t$, that is, to be $C^\infty$
functions in $E$ with values in $\bbR$. A solution to system (\ref{sysdyn}) is a trajectory
curve $\mb{X}(t)$. Since none of components $F_i$ depends explicitely on time, the system is said
to be autonomous. The acceleration vector $\ddot{\mb{X}}$ of a dynamical system can be written as
\begin{equation}
  \label{jaco}
  \ddot{\mb{X}} = {\mb{J}} \dot{\mb{X}}
\end{equation}
where ${\mb{J}}$ is the functional Jacobian matrix of the system.

Trajectory curves integral to dynamical system (\ref{sysdyn}) can be viewed as curves in a $n$-dimensional Euclidean space.
They possess local metric properties, namely {\it curvatures}, which can be analytically
deduced from the so-called Fr\'enet formula (see next section) since only time derivatives of the
trajectory curves are involved in the definition of curvature. For dynamical systems in $\bbR^2$
and $\bbR^3$ the concept of curvature may be exemplified. A curve in $\bbR^2$ is a plane curve
which has a {\it torsion} vanishing identically. A curve in $\bbR^3$ has two curvatures, named
{\it curvature} and {\it torsion}, which are also known as first and second curvatures,
respectively. Curvature measures the curve deviation from a straight line in the neighborhood of
any of its points. Roughly, torsion measures magnitude and sense of the curve deviation from the
osculating plane defined as the plane spanned by the instantaneous velocity and acceleration
vectors. Physically, a straight line can be deformed into any 3D curve by bending (curvature)
and twisting (torsion). A curve in $n$-dimensional Euclidean space ($n>3$) has ($n-1$) curvatures which may be
computed using a Gram-Schmidt procedure.

\newpage

The set of points where the curvature of the flow, that is, the curvature of the trajectory of any
$n$-dimensional dynamical system, vanishes defines a ($n-1$) dimensional invariant manifold. The
flow curvature manifold is thus defined by
\begin{equation}
  \label{phidef}
  \begin{array}{rl}
     \phi(\mb{X}) &
     = \dot{\mb{X}} \cdot
       \left( \displaystyle \ddot{\mb{X}} \wedge \stackrel{...}{\mb{X}} \wedge ... \wedge
          \stackrel{n}{\mb{X}}
       \right)
     = \mbox{det }
    \left( \displaystyle \dot{\mb{X}}, \ddot{\mb{X}},  \stackrel{...}{\mb{X}},  ...  \stackrel{n}{\mb{X}}
       \right) = 0
  \end{array}
\end{equation}
where $\stackrel{n}{\mb{X}}$ represents the time derivatives of $\mb{X}$. For a proof, see
\cite{Gin08}. For a three-dimensional dynamical system, the sets of points where curvature of the
flow vanishes defines a two-dimensional invariant manifold whose analytical equation reads
\begin{equation}
  \phi (\mb{X}) =
  \dot{\mb{X}} \cdot \left( \displaystyle \ddot{\mb{X}} \wedge \stackrel{...}{\mb{X}} \right) =
  \mbox{det }
    \left( \displaystyle \dot{\mb{X}}, \ddot{\mb{X}},  \stackrel{...}{\mb{X}} \right) = 0 \, .
\end{equation}
In this case, the manifold is defined by points where the {\it torsion} vanishes.

\begin{figure}[ht]
  \begin{center}
 \includegraphics[height=4.3cm]{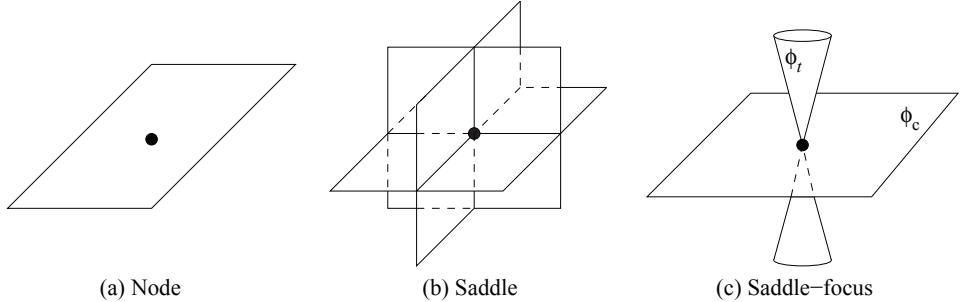} \\[-0.2cm]
    \caption{Generic shapes of the flow curvature manifold in the neighborhood of fixed points.
There is a time dependent component only for a saddle-focus fixed point.}
     \label{torfipo}
  \end{center}
\end{figure}

Differentiating (\ref{jaco}) with respect
to time $t$ leads to
\begin{equation}
  \stackrel{...}{\mb{X}} = {\mb{J}} \ddot{\mb{X}} + \frac{{\rm d}{\mb{J}}}{{\rm d}t} \dot{\mb{X}} \, .
\end{equation}
Inserting this expression into (\ref{phidef}), we obtain
\begin{equation}
  \phi (\mb{X}) =
  \underbrace{\dot{\mb{X}} \cdot \left( {\mb{J}} \dot{\mb{X}} \wedge
            {\mb{J}} \ddot{\mb{X}} \right)}_{\phi_c}
  + \underbrace{\dot{\mb{X}} \cdot \left( \displaystyle \ddot{\mb{X}} \wedge
            \frac{{\rm d}{\mb{J}}}{{\rm d}t} \dot{\mb{X}} \right)}_{\phi_t}
\end{equation}
where $\phi_c$ is the time independent component and $\phi_{t}$ the time dependent component
\cite{Gin08}. Since $\phi_c$ does not contain time derivative of $J$ it is associated with the linear component of the vector field and
$\phi_{t}$ with the nonlinear component. In the neighborhood of fixed points $\mb{X}^*$, the time
independent component of the flow curvature manifold corresponds to the osculating plane
\cite{Gin08}. As a consequence, the attractor takes
the shape of $\phi_c$ in this neighborhood because  the osculating plane cannot be crossed by a
trajectory. This results from the fact that the osculating plane is invariant with respect to the
flow. In all cases, the flow curvature manifold is thus made of a plane parallel to the osculating
plane. In the case of a saddle, time-independent component $\phi_c$ is also made of two additional
transverse planes (Fig.\ \ref{torfipo}b). The fixed point is at the intersection of these three
planes. The two complex conjugated eigenvalues of saddle-focus fixed points induce a non null
time-dependent component which takes the form of two elliptic parabolo\"ids, one associated with the
each branch of the 1D manifold of the fixed point (Fig.\ \ref{torfipo}c). Fixed points of a
saddle-focus type are the only ones with a non-null time-dependent component $\phi_t$.

\section{Rössler-like systems}
\label{Rosso}

The way according which the flow curvature manifold structures the flow is now illustrated for
R\"ossler-like systems, that is, for systems which have R\"ossler-like attractors for their solutions.

\subsection{Systems with two fixed points}

Let us start with the original Rössler system \cite{Ros76a}:
\begin{equation}
  \label{ros76}
  \left\{
    \begin{array}{l}
      \dot{x} = -y - z \\
      \dot{y} = x + ay \\
      \dot{z} = b + z (x-c) \, .
    \end{array}
  \right.
\end{equation}
We choose to center the R\"ossler system but this is not compulsory for our analysis.  The
R\"ossler system is thus centered through a rigid displacement, that is, the inner fixed point,
$F_-$, is moved to the origin of the phase space $\bbR^3 (x,y,z)$. In the translated coordinate
system, the equations for the centered system are
\begin{equation}
\label{ros79cen}
  \left\{
    \begin{array}{l}
      \dot{x} = -y - z -y_- -z_- \\
      \dot{y} = x + ay + x_- + a y_- \\
      \dot{z} = b + z (x+x_- -c) + z_- x + z_- (x-c)
    \end{array}
  \right.
\end{equation}
where $\frac{x_-}{a} = -y_- = z_- = \frac{c-\sqrt{c^2-4ab}}{2a}$
are the coordinates of the inner fixed
point of the Rössler system (\ref{ros76}). The system may then be rewritten as:
\begin{equation}
  \label{roscent}
  \left\{
    \begin{array}{l}
      \dot{x} = -y - z \\
      \dot{y} = x + ay \\
      \dot{z} = \tilde{b}x + z (x-\tilde{c}) \, .
    \end{array}
  \right.
\end{equation}
where $\tilde{b} = z_-$ and $\tilde{c} = c-x_-$. This centered Rössler system has one fixed point
$F_-$ located at the origin of the phase space and another one located at
\begin{equation}
  F_+ =
  \left|
    \begin{array}{l}
      x_+ = \tilde{c} - a \tilde{b} \\[0.2cm]
      \displaystyle
      y_+ = - \frac{x_+}{a} \\[0.2cm]
      \displaystyle
      z_+ = \frac{x_+}{a}
    \end{array}
  \right. \, .
\end{equation}
The structure of the flow near the origin and along the $x$-$y$ plane is governed to a large
extent by the unstable fixed point $F_-$ (previously designated as the inner fixed point). This causes
the flow to ``spiral around'' this point.  On a larger scale, the flow in the R\"ossler attractor
wraps around the one-dimensional unstable manifold associated with the outer fixed point $F_+$.

\begin{figure}[ht]
  \begin{center}
    \begin{tabular}{ccc}
      \includegraphics[height=5.0cm]{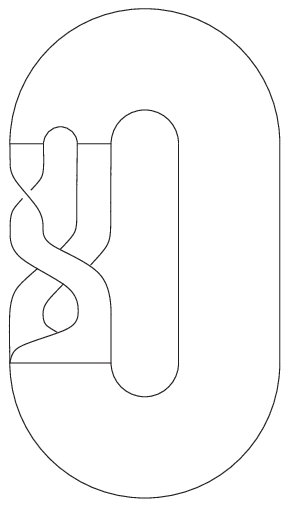} & ~~~ &
      \includegraphics[height=5.0cm]{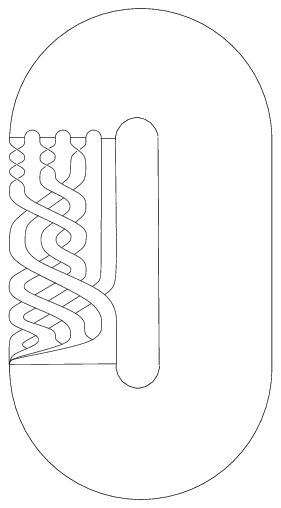} \\
       (a) Two branch template & & (b) Four branch template \\[-0.2cm]
    \end{tabular}
    \caption{Templates for two different chaotic attractors solution to the Rössler system.
Typical parameter values: $b=2$ and $c=4$. Template (a) is obtained with $a=0.432$ and template
(b) for $a=0.52$.}
   \label{templates}
  \end{center}
\end{figure}

\newpage

\begin{figure}[ht]
  \begin{center}
    \begin{tabular}{ccc}
      \includegraphics[height=5.0cm]{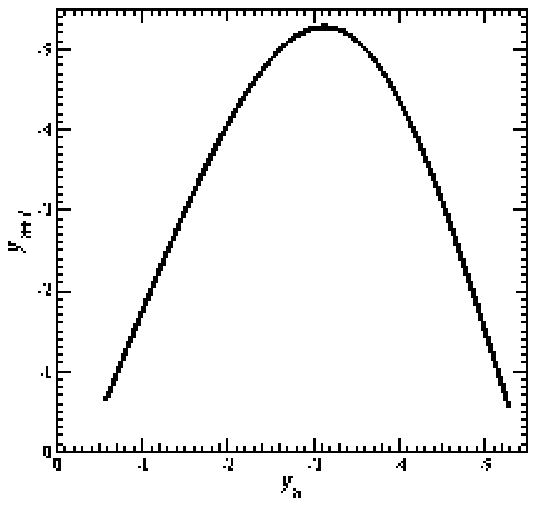} & ~~~ &
      \includegraphics[height=5.0cm]{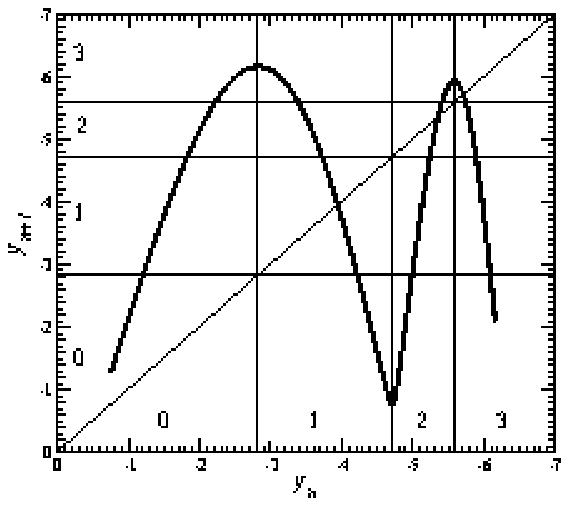} \\[-0.3cm]
      (a) $a=0.432$ & & (b) $a=0.52$ \\[-0.2cm]
    \end{tabular}
    \caption{First-return map to a Poincaré section of two different chaotic attractors solution to
the Rössler system. Typical parameter values: $b=2$ and $c=4$.}
   \label{mapros}
  \end{center}
\end{figure}

The simplest chaotic attractor solution to the R\"ossler system has a topology which can be
described by a template with two branches as shown in Fig. \ref{templates}a \cite{Let95a}. Its
first-return map to a Poincaré section presents two monotonic branches (Fig.\ \ref{mapros}a).
When parameter $a$ is increased, the attractor after a sequence of bifurcations becomes of funnel
type, that is, characterized by a first-return map to a Poincar\'e section with many monotone
branches (four in the case shown in Fig.\ \ref{mapros}b). The template has therefore two
additional branches (Fig.\ \ref{templates}b) compared to the previous template (Fig.
\ref{templates}a). In order to describe the way in which monotonic branches are developed and
visited, a partition of the attractor can be defined according to the critical points (extrema)
of the first-return map (Fig.\ \ref{mapros}b). A transition matrix is thus defined according
to the panels where at least one point can be found. In the case of the first-return map shown in
Fig.\ \ref{mapros}b, all panels are visited and the corresponding transition matrix is
\begin{equation}
  \Gamma =
  \left[
    \begin{array}{cccc}
      ~1~ & ~1~ & ~1~ & ~1~ \\
      1 & 1 & 1 & 1 \\
      1 & 1 & 1 & 1 \\
      1 & 1 & 1 & 1
    \end{array}
  \right] \, .
\end{equation}
A detailed study of the R\"ossler attractor can be found in \cite{Let95a}.

\begin{figure}[ht]
  \begin{center}
    \includegraphics[height=6.0cm]{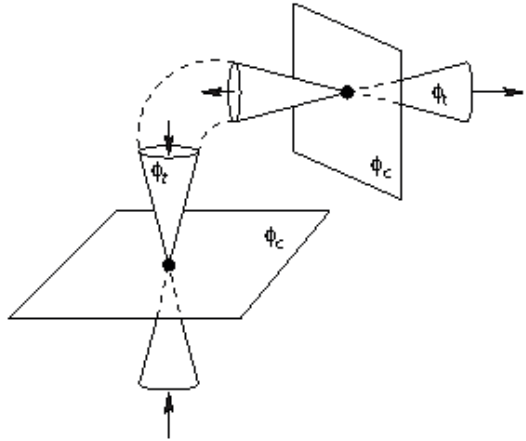} \\[-0.2cm]
     \caption{Scheme of the flow curvature manifold for the Rössler attractor. The two elliptic
parabolo\"ids from the fixed points are joined to form a single closed ellipso\"id.}
    \label{nutoros}
  \end{center}
\end{figure}

According to the generic shapes for the time-independent component of the flow curvature
manifold identified in the previous section, a scheme of the flow curvature manifold can be drawn
as shown in Fig.\ \ref{nutoros}. The inner fixed point $F_-$ has a plane associated with its
unstable 2D manifold and an elliptic parabolo\"id centered on its stable 1D manifold. The outer fixed
point $F_+$ has a elliptic parabolo\"id associated with its unstable 1D manifold and a plane
corresponding to the
stable 2D manifold. In all systems investigated in this paper the inner fixed point has a 2D
unstable manifold and those associated with the outer fixed point is 1D. To our knowledge, there is
no continuous dynamical system producing an attractor topologically equivalent to the Rössler
attractor, and surrounding a fixed point with a 2D stable manifold.

\begin{figure}[ht]
  \begin{center}
    \begin{tabular}{cc}
      \includegraphics[height=5.4cm]{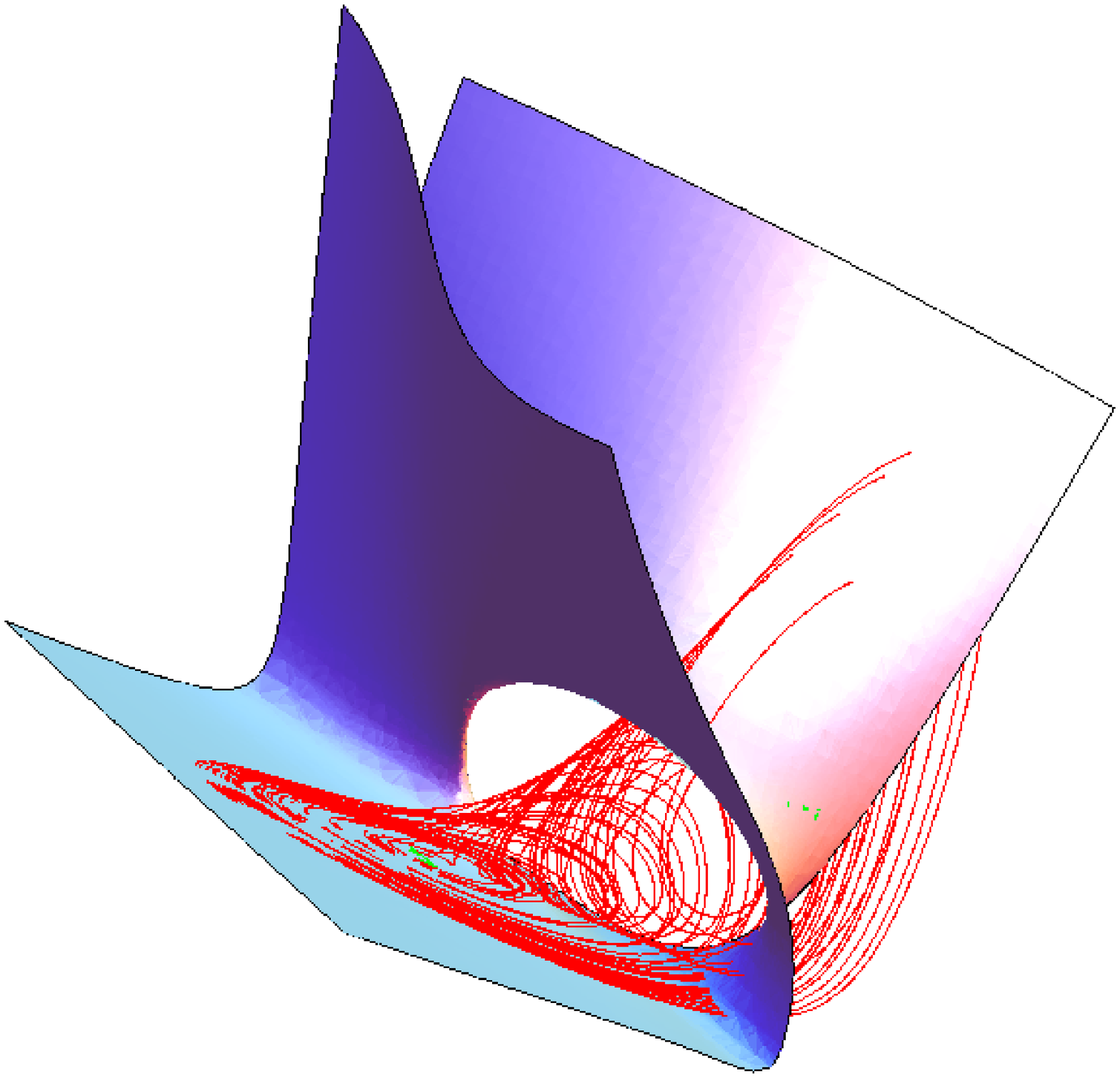} &
      \includegraphics[height=4.4cm]{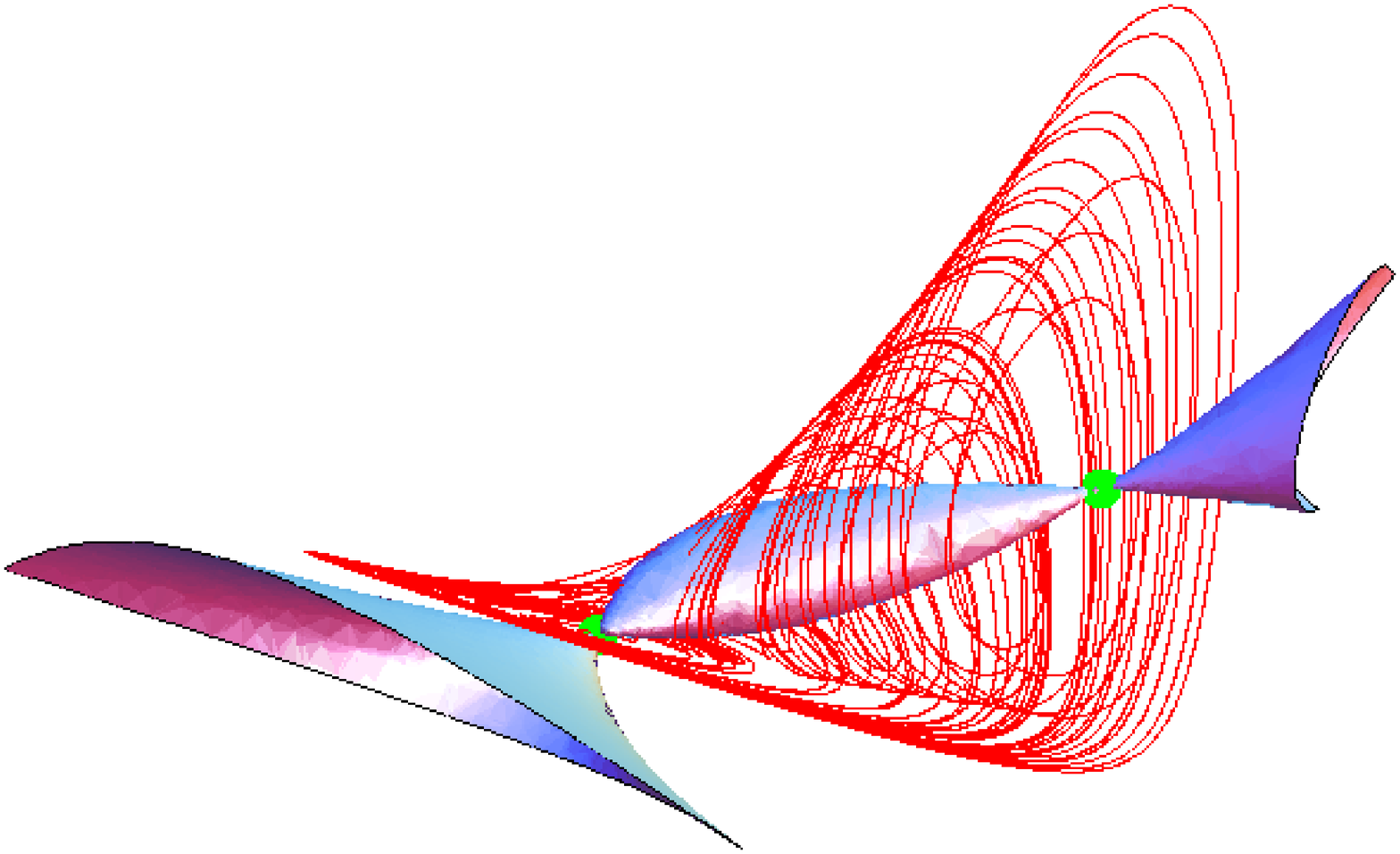} \\[-0.0cm]
      (a) Time-independent component $\phi_c$ &
      (b) Time-dependent component $\phi_t$ \\[-0.2cm]
    \end{tabular}
    \caption{The two components of the flow curvature manifold $\phi$ for the Rössler system
with parameter values: $a=0.556$, $b=2$ and $c=4$.}
   \label{rosunuto}
  \end{center}
\end{figure}

The two components of the flow curvature manifold of the Rössler system are shown in Fig.\
\ref{rosunuto}. As expected, in the neighborhood of the inner fixed point, time dependent
component of the flow curvature manifold is tangent to the osculating plane, that is, nearly
parallel to the $x$-$y$ plane. Component $\phi_t$ presents an elliptic parabolo\"id at each side of
the 2D manifolds of the fixed points. Between the two fixed points, these elliptic parabolo\"ids are
joined to form
a closed ellipso\"id (Fig.\ \ref{rosunuto}b). The trajectory wraps around a significant part of
this closed ellipso\"id. Close to the inner fixed point, the trajectory
crosses component $\phi_t$. Note that the boundary of the non visited neighborhood of the inner
fixed point roughly corresponds to the location where the trajectory crosses component $\phi_t$.
Such an intersection between the trajectory and component $\phi_t$ could be an explaination to the
limitation to the development of the dynamics. According to such an assumption, such a crossing
could be responsible for the pruning of periodic orbits observed in the neighborhood of the inner
fixed point \cite{Let95a}. This is confirmed by the fact that, for $a=0.43295$, the trajectory
visits the neighborhood of the inner fixed point and does not intersect component $\phi_t$.

\begin{figure}[ht]
  \begin{center}
      \includegraphics[height=5.4cm]{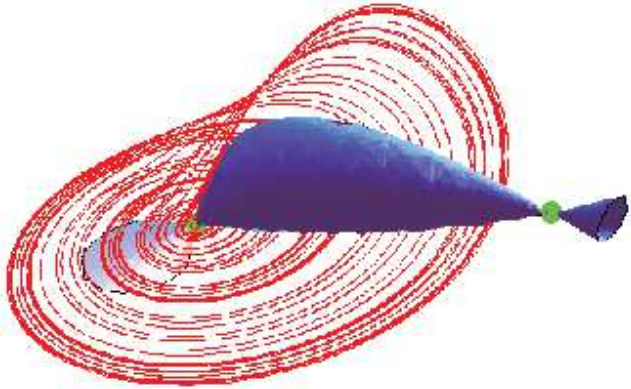} \\[-0.2cm]
      \caption{Time-dependent component $\phi_t$ of the flow curvature manifold for the
Rössler system with parameter values: $a=0.43295$, $b=2$ and $c=4$.}
    \label{rosflocur}
  \end{center}
\end{figure}

\begin{figure}[htbp]
  \begin{center}
    \begin{tabular}{cc}
      Chaotic attractor & First-return map \\[-0.0cm]
      \includegraphics[height=6.4cm]{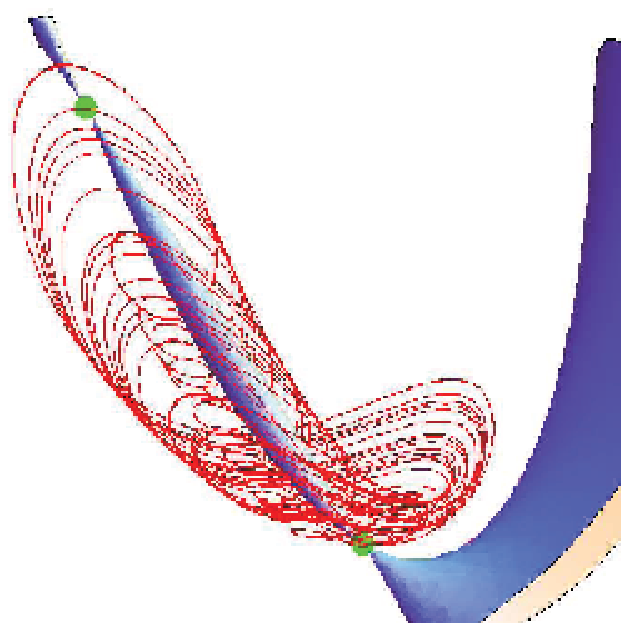} &
      \includegraphics[height=6.4cm]{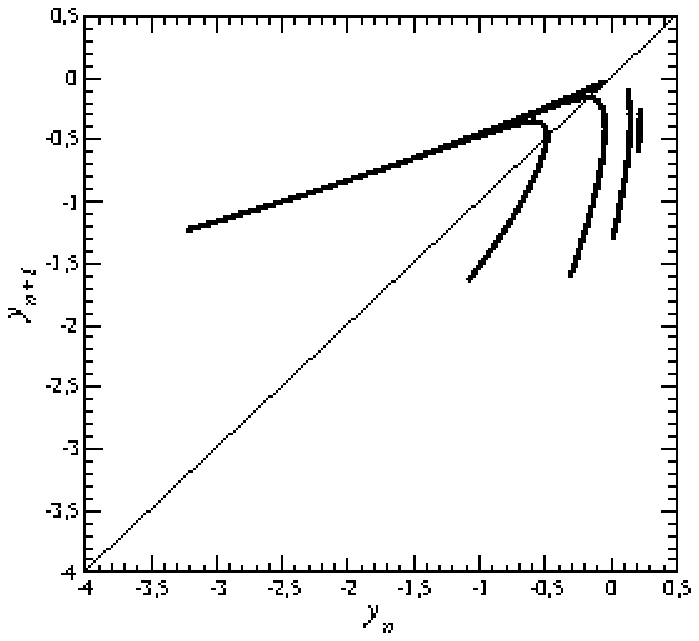} \\[-0.2cm]
      \multicolumn{2}{c}{(a) Sprott system F, $a=0.5$ } \\[0.2cm]
      \includegraphics[height=6.4cm]{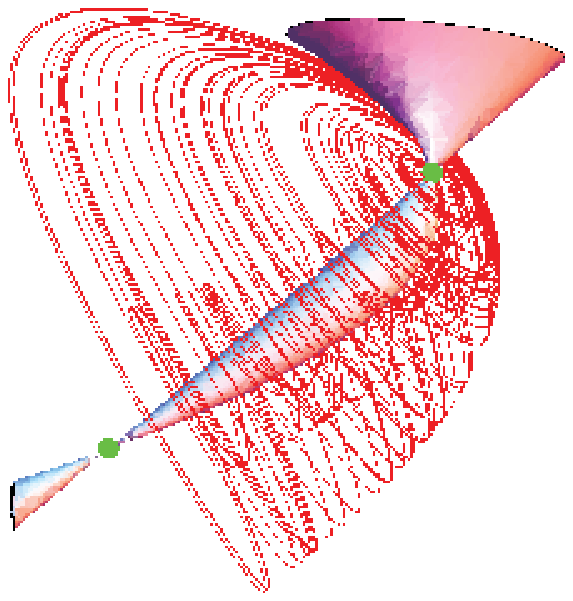} &
      \includegraphics[height=6.4cm]{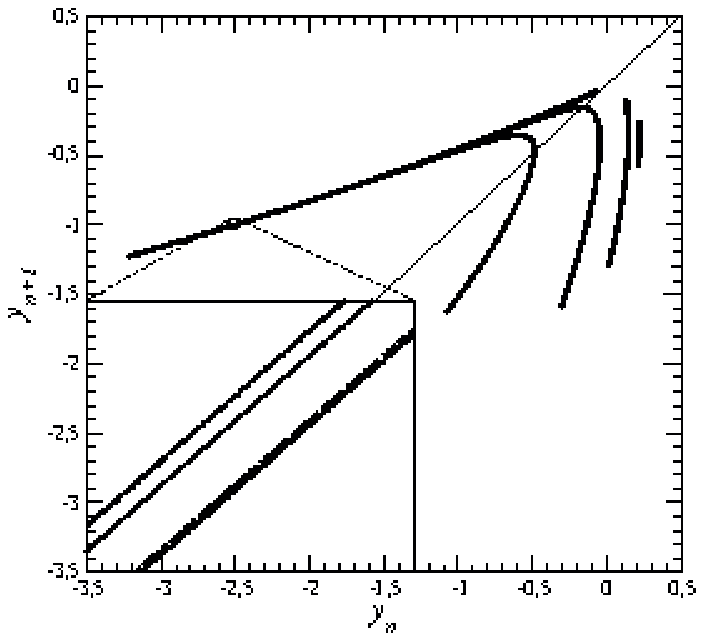} \\[-0.2cm]
      \multicolumn{2}{c}{(b) Sprott system H, $a=0.5$} \\[0.2cm]
      \includegraphics[height=6.4cm]{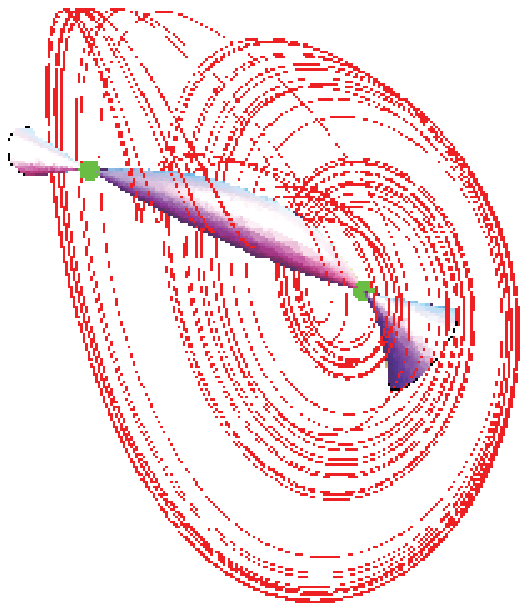} &
      \includegraphics[height=6.4cm]{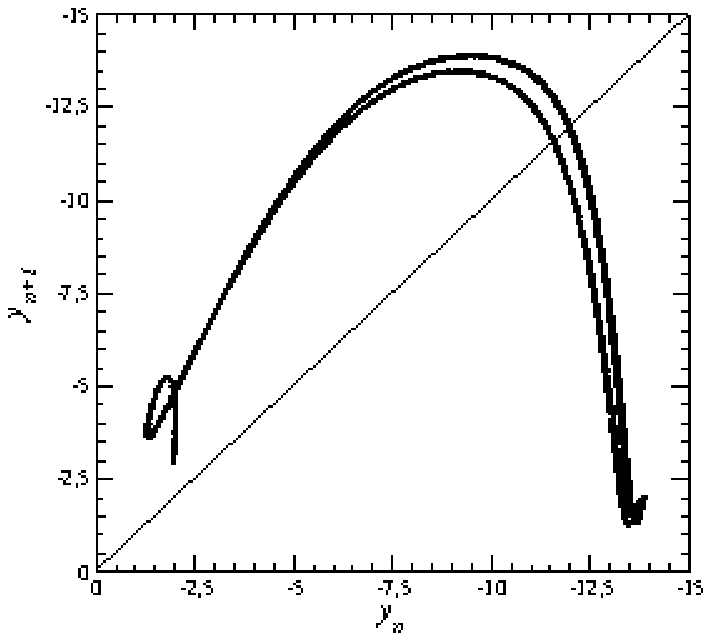} \\[-0.2cm]
      \multicolumn{2}{c}{(c) Sprott system Q, $a=3.1$ and $b=0.5$ } \\[-0.2cm]
    \end{tabular}
    \caption{Chaotic solutions to Sprott system F, H and Q.}
    \label{sprottF}
  \end{center}
\end{figure}

Nine other Rössler-like systems were investigated. All the other systems investigated in the
subsequent part of this paper can be written under the general form

\begin{equation}
  \label{genform}
  \left\{
    \begin{array}{l}
      \dot{x} = a_2 y + a_3 z + a_4 xz + a_5 z^2 \\[0.1cm]
      \dot{y} = b_1 x + b_2 y + b_3 z + b_4 y^2 + b_5 z^2 \\[0.1cm]
      \dot{z} = c_1 x + c_2 y + c_3 z + c_4 xy + c_5 xz + c_6 x^2 \\[0.1cm]
    \end{array} \, .
  \right.
\end{equation}
Only coefficents $a_i$, $b_j$ and $c_k$ are reported in Tab.\ \ref{syscoe}. In all of these systems
but one, the elliptic parabolo\"ids emerging from the fixed points form a closed ellipso\"id
(Figs \ref{rosunuto}b and \ref{sprottF}).

{\small
\begin{table}[ht]
  \begin{center}
    \caption{Specific coefficients of each system here investigated. Compared to their original
form as published in \cite{Ros76a} and \cite{Spr94}, each system was centered, that is, the
inner fixed point was located at the origin of the phase space.}
    \label{syscoe}
    \begin{tabular}{llcccccccccccccccccccccc}
      \\[-0.3cm]
      \hline \hline
      \\[-0.3cm]
      & & \multicolumn{4}{c}{$\dot{x}=$} &
      \multicolumn{5}{c}{$\dot{y}=$} &
      \multicolumn{6}{c}{$\dot{z}=$} \\
      & & $y$ & $z$ & $xz$ & $z^2$ & $x$ & $y$ & $z$ & $y^2$ & $z^2$ &
$x$ & $y$ & $z$ & $xy$ & $xz$ & $x^2$ \\
      System & &
    $a_2$ & $a_3$ & $a_4$ & $a_5$ &
    $b_1$ & $b_2$ & $b_3$ & $b_4$ & $b_5$ &
    $c_1$ & $c_2$ & $c_3$ & $c_4$ & $c_5$ & $c_6$ \\[0.1cm] \hline
      \\[-0.3cm]
Rössler & & -1 & -1 & 0 & 0 & +1 & $+a$ & 0 & 0 & 0 & $\tilde{b}$ & 0 & $-\tilde{c}$
& 0 & +1 & 0 \\
Sprott F & & -1 & +1 & 0 & 0 & +1 & $+a$ & 0 & 0 & 0 & 0 & 0 & -1 & 0 & 0 & +1 \\
Sprott G & & -1 & +1 & 0 & 0 & +1 & $+a$ & 0 & 0 & 0 & 0 & 0 & $-b$ & +1 & 0 & 0 \\
Sprott H & & -1 & 0 & 0 & +1 & +1 & $+a$ & 0 & 0 & 0 & +1 & 0 & $-1$ & 0 & 0 & 0 \\
Sprott K & & -1 & 0 & +1 & 0 & +1 & $+a$ & 0 & 0 & 0 & +1 & 0 & $-b$ & 0 & 0 & 0 \\
Sprott M & & -1 & 0 & 0 & 0 & $+a$ & 0 & +1 & 0 & 0 & $+b$ & 0 & -1 & 0 & 0 & -1 \\
Sprott O & & +1 & 0 & 0 & 0 & +1 & 0 & -1 & 0 & 0 & +1 & $+a$ & 0 & 0 & +1 & 0 \\
Sprott P & & $+a$ & +1 & 0 & 0 & -1 & 0 & 0 & +1 & 0 & +1 & +1 & 0 & 0 & 0 & 0 \\
Sprott Q & & -1 & 0 & 0 & 0 & $+a$ & $+b$ & 0 & 0 & +1 & +1 & 0 & -1 & 0 & 0 & 0 \\
Sprott S & & +1 & 0 & 0 & 0 & 0 & $-a$ & $-b$ & 0 & 0 & +2 & +1 & 0
& 0 & 0 & +1 \\
      \\[-0.3cm]
     \hline \hline
    \end{tabular}
  \end{center}
\end{table}
}

For instance, we observe that Sprott systems F and H produce well developed
and similar funnel attractors (Figs.\ \ref{sprottF}a and \ref{sprottF}b). For these two systems,
the trajectory wraps around component $\phi_t$ --- and therefore does not cross it --- almost
everywhere between the two fixed points. First-return map to a Poincaré section of attractors
solution to Sprott systems F and H have unusual shapes. Four decreasing monotonous branches are
clearly distinguished and a blow up shows four increasing branches (Figs.\ \ref{sprottF}a and
\ref{sprottF}b). The map has thus eight branches. Such a feature results from the numerical
difficulties in computing a proper Poincaré section. Sprott system Q also does not present a trajectory
crossing component $\phi_t$ too but its funnel structure is less developed (Fig.\ \ref{sprottF}c)
than the one of Sprott systems F and H. In particular, the first-return map has only two branches
(Fig.\ \ref{sprottF}c). The main departure between these systems could be how fast the trajectory
wraps around component $\phi_t$.

In order to roughly quantify this dynamical property, we compute a wrapping number
defined as
\begin{equation}
  W = \frac{\omega}{\lambda_3} D_{F_+ - F_-}
\end{equation}
where $\omega$ is the imaginary part of the complex conjugated eigenvalues of the outer fixed
point, $\lambda_3$ its real eigenvalue and $D_{F_+ - F_-}$ the distance between the two
fixed points $F_+$ and $F_-$. For the three Sprott systems F, H and Q, we obtained
$W_F = 59.4$, $W_H = 48.5$ and $W_Q = 0.2$, respectively. Obviously, the trajectory solution
to Sprott system Q wraps more slowly than trajectories solution to Sprott systems F and H.
The dynamics of Sprott system Q is therefore less developed. In this case, such a limitation results
from the eigenvalues of the outer fixed point.

\begin{figure}[ht]
  \begin{center}
    \begin{tabular}{cc}
      \includegraphics[height=6.0cm]{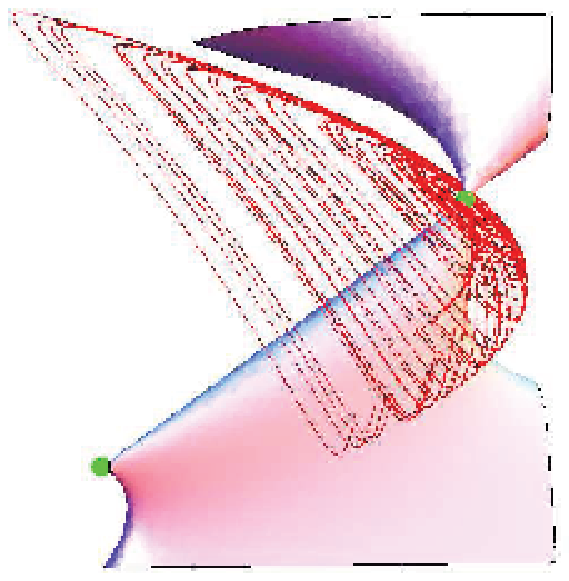} &
      \includegraphics[height=6.0cm]{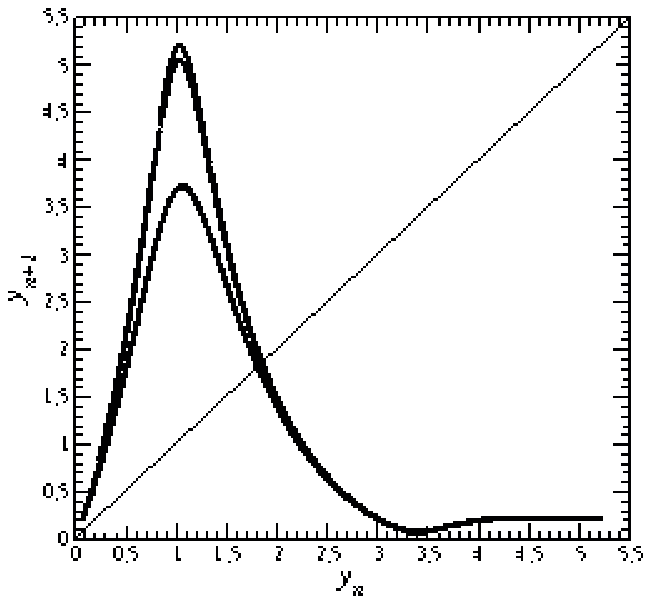} \\[-0.2cm]
      (a) Chaotic attractor & (b) First-return map \\[-0.2cm]
    \end{tabular}
    \caption{Chaotic behavior solution to Sprott system K. Parameter values: $a=0.35$ and $b=0.5$.
Component $\phi_t$ is not a closed ellipso\"id due to a singularity which appears when solving
$\phi (x,y,z)=0$. By applying the implicit function theorem, we can express $z=\Psi (x,y)$
in terms of $x$ and $y$ where there is a singularity in $x$ inducing numerical artifacts.}
    \label{sprottK}
  \end{center}
\end{figure}

\begin{figure}[ht]
  \begin{center}
    \begin{tabular}{cc}
      \includegraphics[height=5.0cm]{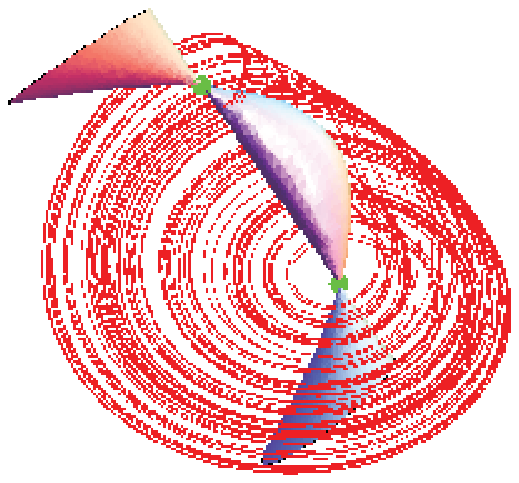} &
      \includegraphics[height=6.0cm]{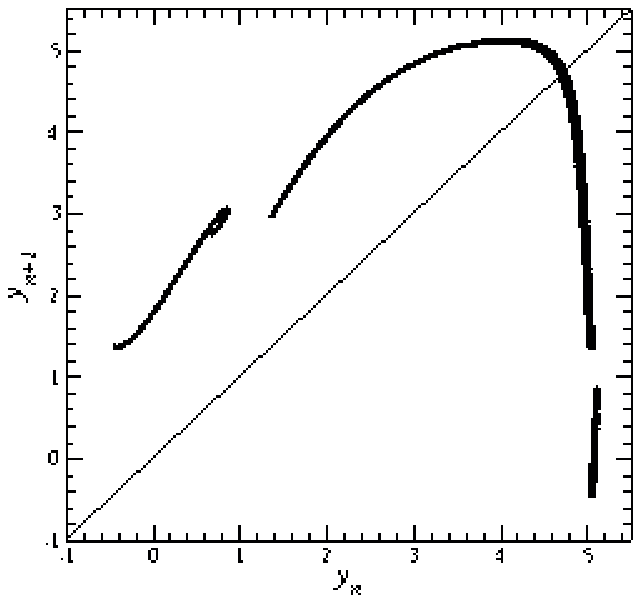} \\[-0.2cm]
      (a) Chaotic attractor & (b) First-return map \\[-0.2cm]
    \end{tabular}
    \caption{Chaotic behavior solution to Sprott system S. Parameter values: $a=0.99$ and $b=3.8$.}
    \label{sprottS}
  \end{center}
\end{figure}

It must be pointed out that the eigenvalues of the outer fixed point do not explain the development
of all attractors investigated here. Indeed, when wrapping numbers $W$ are computed for the five
other Sprott systems reported in Tab.\ \ref{syscoe}, we got
\[
  \begin{array}{cc}
    W_S = 3.8 < W_O = 4.3 < W_P = 8.5 < W_M = 14 < ... \\[0.1cm]
   ... < W_G = 21.3 < W_K = 27.1 \, .
  \end{array}
\]
In particular, $W_K$ is significantly greater than $W_S$ but the attractor solution to Sprott
system S has an attractor (Fig.\ \ref{sprottS}) which is not significantly more developed than the
attractor solution to Sprott system K (Fig. \ref{sprottK}): the latter presents a unimodal map
(Fig.\ \ref{sprottK}b) and the former a three branches map (Fig.\ \ref{sprottS}b) where the
third branch is rather small. Moreover, $W_K$ is around the
half of $W_F$ and a more developed dynamics (at least four branches) was expected. The major
ingredient, observed in Sprott systems K and S but not in systems F and H, is that the trajectory
intersects the time-dependent component $\phi_t$. Such an intersection is viewed as being the main
reason for the limitation of the dynamics, that is, of the number of monotonous branches in the
first-return map.

\begin{figure}[htbp]
  \begin{center}
    \begin{tabular}{cc}
      \includegraphics[height=5.0cm]{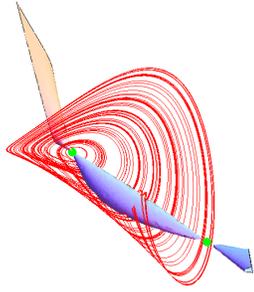} &
      \includegraphics[height=5.0cm]{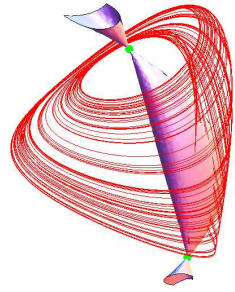} \\[-0.0cm]
      (a) Sprott system G & (b) Sprott system M \\
      $a=0.42$ and $b=1.29$ & $a=1.95$ and $b=1.65$ \\[0.2cm]
      \includegraphics[height=5.5cm]{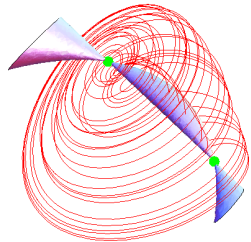} &
      \includegraphics[height=5.0cm]{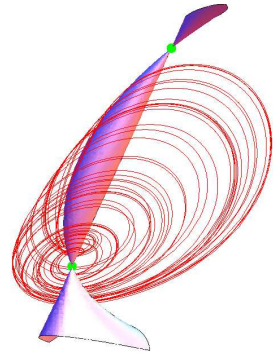} \\
      (c) Sprott system O & (d) Sprott system P \\
      $a=2.67$ and $b=0.5$ & $a=2.68$  \\[-0.2cm]
    \end{tabular}
    \caption{Chaotic attractors solution to Rössler-like systems with their two fixed points
(designated by green circles in the figures) and the time-dependent component of their flow
curvature manifold. Parameter values correspond to the most developed attractor we identified for
each system.}
   \label{sprottG}
  \end{center}
\end{figure}

The structure of Rössler-like attractors therefore depends on the fixed
points (and their eigenvalues) and, the interplay between the flow curvature manifold and the
trajectory. The core of the time-dependent component $\phi_t$ can be considered as an
axis around which the trajectory wraps when there is no intersection between the trajectory and
component $\phi_t$.  The four remaining Sprott systems with two fixed points are quite similar to
the case of Sprott system S (Fig.\ \ref{sprottG}).

\begin{figure}[ht]
  \begin{center}
    \begin{tabular}{ccc}
      \includegraphics[height=2.5cm]{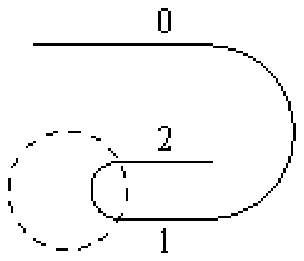} & ~~ &
      \includegraphics[height=3.5cm]{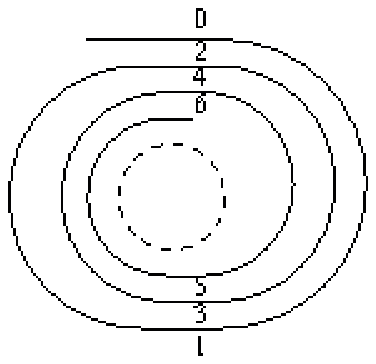} \\
      (a) With intersection & & (b) Without intersection \\[-0.2cm]
    \end{tabular}
    \caption{Scheme of the transverse structure to the flow observed in a Poincaré section with
component $\phi_t$ (dashed line). When the trajectory wraps around $\phi_t$ the number of
branches in the first-return map --- or equivalently in the template --- is limited by the ratio
$W$ (b). Once there is an intersection between the trajectory and component $\phi_t$, the number
of branches can no longer increase because the trajectory is not longer in the neighborhood of
component $\phi_t$ (a).}
    \label{withint}
  \end{center}
\end{figure}

In fact, when the trajectory intersects component $\phi_t$, it presents a folding rather than
a wrapping structure. Once the
trajectory crossed component $\phi_t$ and described a fold, it is no longer located in a zone of
the phase space where there is a structure (component $\phi_t$) around which it can wrap
(Fig.\ \ref{withint}a). The corresponding attractor can no longer develop new branches and the
``funnel'' type is quite limited (most often three branches in the first-return map). The
probability having an intersection between the trajectory and component $\phi_t$ seems to be
greater than not.  This would explain why limited funnel attractors are more often observed.

\subsection{Systems with a single fixed point}
\label{others}

In his exhaustive search procedure, Sprott also found systems with a single fixed point. Seven of
them will be investigated in this section. Once these systems were centered, they have the general
form (\ref{genform}) and their c\oe fficients are reported in Tab.\ \ref{syscod}. One system
proposed by Thomas \cite{Tho07} and two by Malasoma \cite{Mal02} were also considered. For all of
these systems, parameter values used for this study correspond to the most developed chaotic
attractor we observed in these systems.

\begin{table}[ht]
  \begin{center}
    \caption{Specific coefficients for the systems with a single fixed point here investigated.
The last four systems --- Sprott systems L and N, and Malasoma systems A and B --- produce an
inverted Rössler-like chaotic attractor.}
    \label{syscod}
    \begin{tabular}{lccccccccccccccccccccccc}
      \\[-0.2cm]
      \hline \hline
      \\[-0.2cm]
      & & \multicolumn{2}{c}{$\dot{x}=$} &
      \multicolumn{4}{c}{$\dot{y}=$} &
      \multicolumn{5}{c}{$\dot{z}=$} \\
      & & ~$y$~ & ~$z$~ & ~$x$~ & ~$y$~ & ~$z$~ & ~$z^2$~ &
~$x$~ & ~$y$~ & ~$z$~ & ~$xy$~ & ~$y^2$~ \\
      System & &
    $a_2$ & $a_3$ &
    $b_1$ & $b_2$ & $b_3$ & $b_5$ &
    $c_1$ & $c_2$ & $c_3$ & $c_4$ & $c_7$ \\[0.1cm] \hline
      \\[-0.3cm]
Sprott D & & -1 & 0 & +1 & 0 & +1 & 0 & 0 & +1 & $+a$ & 0 & 1 \\
Sprott I & & $-a$ & 0 & +1 & 0 & +1 & 0 & +1 & 0 & -1 & 0 & +1 \\
Sprott J & & $+a$ & 0 & -1 & 0 & +1 & 0 & +1 & +1 & $-a$ & 0 & 0 \\
Sprott R & & -1 & 0 & 0 & 0 & +1 & 0 & $+a$ & $-\frac{b}{a}$ & -1 & +1 & 0 \\
Thomas & & +1 & 0 & -1 & $+a$ & -1 & 0 & 0 & 0 & $-c$ & 0 & +1 \\[0.1cm] \hline
Sprott L & & -1 & 0 & $+a$ & 0 & +1 & 0 & 0 & $+2b$ & -1 & 0 & $+b$ \\
Sprott N & & $-a$ & 0 & +1 & 0 & $+\frac{2}{a}$ & +1 & 0 & +1 & $-a$ & 0 & 0 \\
Malasoma A & & +1 & 0 & 0 & $-a$ & +1 & 0 & -1 & 0 & 0 & +1 & 0 \\
Malasoma B & & 0 & +1 & 0 & $-a$ & +1 & 0 & -1 & 0 & 0 & +1 & 0 \\[0.1cm]
      \\[-0.3cm]
     \hline \hline
    \end{tabular}
  \end{center}
\end{table}

\begin{figure}[ht]
  \begin{center}
    \begin{tabular}{cc}
      \includegraphics[height=5.0cm]{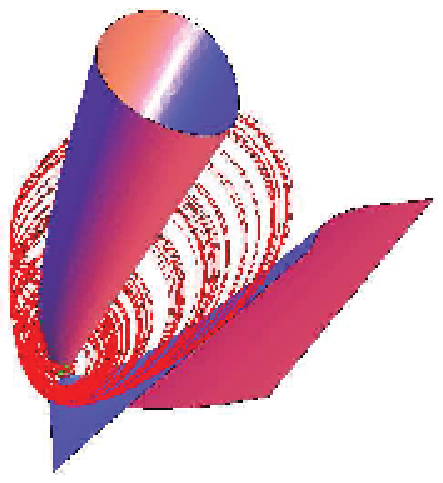} &
      \includegraphics[height=5.0cm]{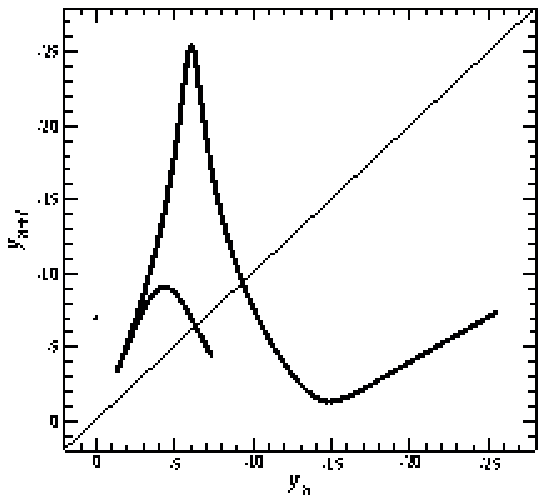} \\[-0.2cm]
      \multicolumn{2}{c}{(a) Sprott system J, $a=1.76$}  \\[0.2cm]
      \includegraphics[height=5.0cm]{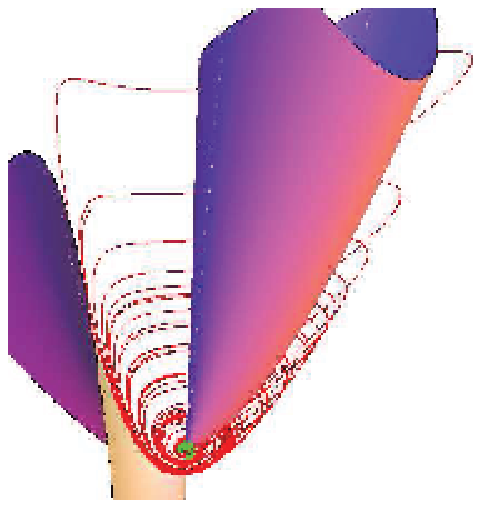} &
      \includegraphics[height=5.0cm]{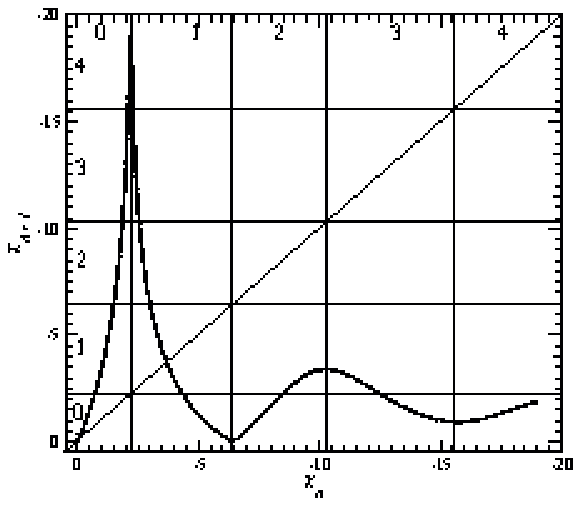} \\[-0.2cm]
      \multicolumn{2}{c}{(b) Thomas system, $a=0.28$ and $c=2$} \\[-0.2cm]
    \end{tabular}
    \caption{Two systems with a single fixed point. The trajectory crosses component $\phi_t$ of
the flow curvature manifold. This limits the development of the attractor.}
    \label{sprottJ}
  \end{center}
\end{figure}

The Sprott system J and Thomas system present a time-dependent component which is crossed by the
trajectory (Fig.\ \ref{sprottJ}). Their attractors are therefore not so developed as in the
previous case. As observed for systems with two fixed points, once the trajectory crosses
component $\phi_t$ of the flow curvature manifold, it is no longer possible to continue
to develop the wrapping process. The resulting attractor is slightly more developed than a unimodal
attractor. Sprott system J presents five branches, two of them being under the first two branches
(Fig.\ \ref{sprottJ}a).
In particular, the small increasing branch is quite difficult to distinguish from the first large
increasing branch due to the difficulty of computing a well-defined Poincaré section. Thomas system
is quite similar to Sprott system J. The advantage of Thomas system is that a safe Poincaré section
can be easily computed. As a consequence, its first-return map clearly presents five monotonous
branches (Fig.\ \ref{sprottJ}b). In both cases, there are two well developed branches and three
others that are not very developed. The corresponding transition matrix
\begin{equation}
  \Gamma =
  \left[
    \begin{array}{ccccc}
      ~1~ & ~1~ & ~1~ & ~1~ & ~1~ \\
      1 & 1 & 1 & 1 & 1 \\
      1 & 1 & 0 & 0 & 0 \\
      1 & 1 & 0 & 0 & 0 \\
      1 & 1 & 0 & 0 & 0
    \end{array}
  \right]
\end{equation}
reveals that, for instance, points in branches labelled 2, 3 and 4 are necessarily followed by
points located in the first two branches (labelled 0 and 1, respectively). According to our views,
this feature results from intersections between the trajectory and component $\phi_t$.

\begin{figure}[ht]
  \begin{center}
    \begin{tabular}{cc}
      \includegraphics[height=6.0cm]{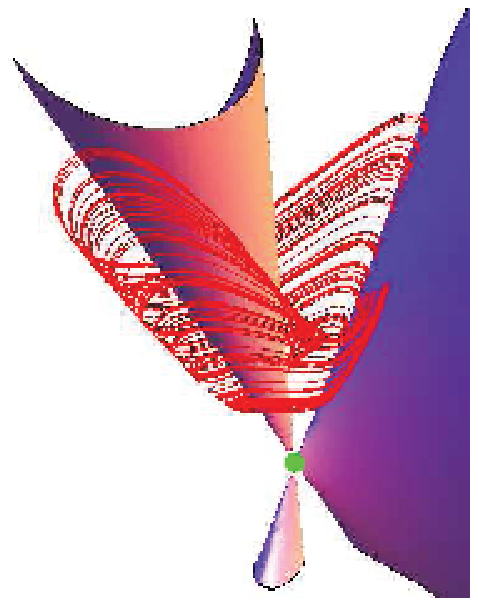} &
      \includegraphics[height=6.0cm]{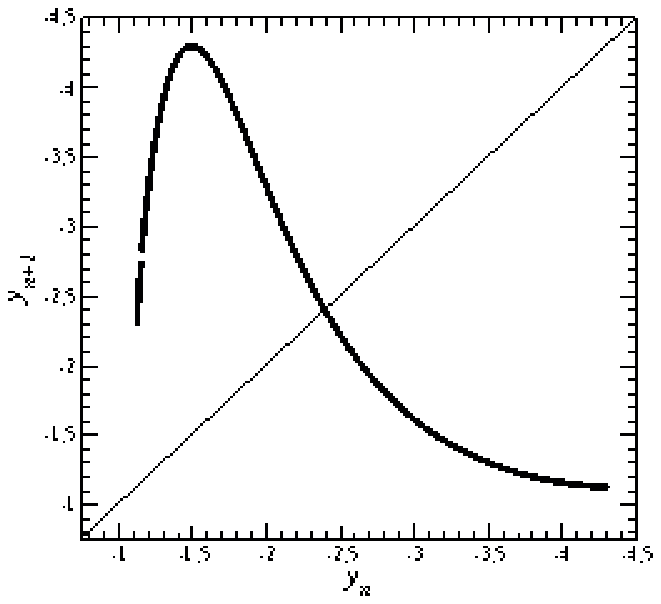} \\[-0.2cm]
      \multicolumn{2}{c}{(a) Sprott system D, $a=2.3$} \\[0.2cm]
      \includegraphics[height=6.0cm]{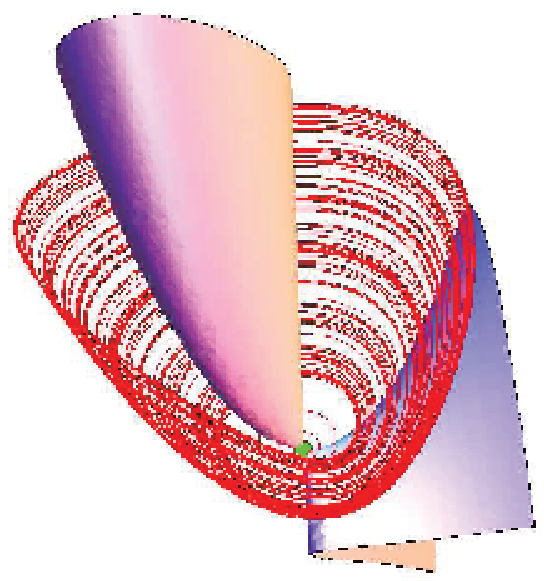} &
      \includegraphics[height=6.0cm]{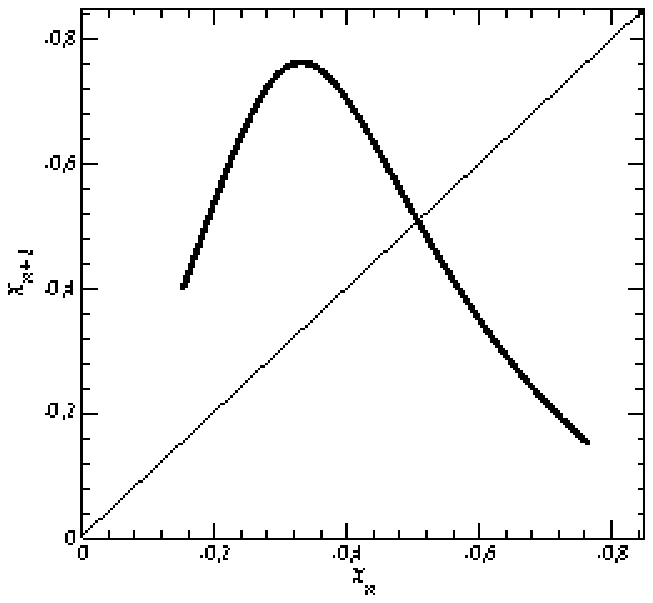} \\[-0.2cm]
      \multicolumn{2}{c}{(b) Sprott system I, $a=0.25$} \\[-0.2cm]
    \end{tabular}
    \caption{Two systems with a single fixed point producing a quite limited chaotic attractor.
System D has two pure imaginary eigenvalues and system I has two complex conjugated eigenvalues
with small real parts (Re$(\lambda_\pm) \approx 0.07$).}
    \label{sprottD}
  \end{center}
\end{figure}

Sprott systems D and I present a different configuration. The trajectory does not intersect
component $\phi_t$ around which it wraps. In the case of Sprott system D (Fig.\ \ref{sprottD}a),
there are numerical artifacts in computing component $\phi_t$ due to a singularity occuring when
solving $\phi (x,y,z)=0$. As a consequence, a spurious part is obtained in addition to the two
elliptic parabolo\"ids usually found. The trajectory intersects the spurious part of component
$\phi_t$ and we can consider that there is no intersection between the trajectory and
component $\phi_t$. What limits the dynamics is in fact the two pure imaginary eigenvalues of the
fixed point which forbid the trajectory visiting the neighborhood of the fixed point. A similar
conclusion is obtained for Sprott system I (Fig.\ \ref{sprottD}b) where the fixed point has two
complex conjugated eigenvalues with very small real parts. In both cases, the development of
the attractors can be understood using fixed point eigenvalues combination with component $\phi_t$.

First-return maps to a Poincaré section of these two attractors present two
monotonic branches that are not fully developed. These two chaotic regimes are therefore less developed
than previous cases that have three monotonic branches in their first-return maps.

\begin{figure}[ht]
  \begin{center}
    \begin{tabular}{cc}
      \includegraphics[height=5.5cm]{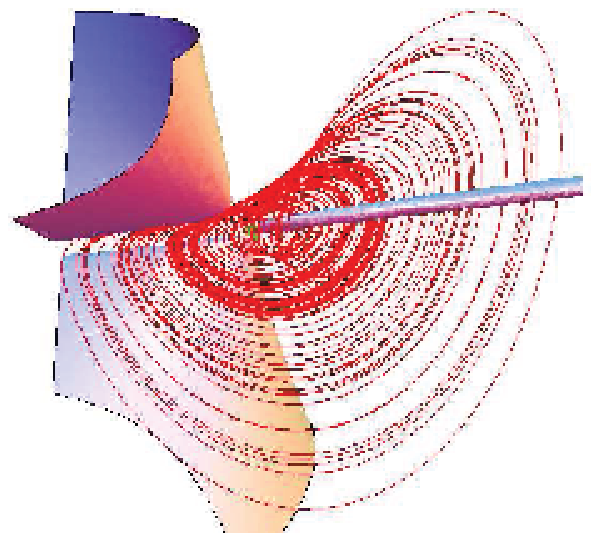} &
      \includegraphics[height=5.5cm]{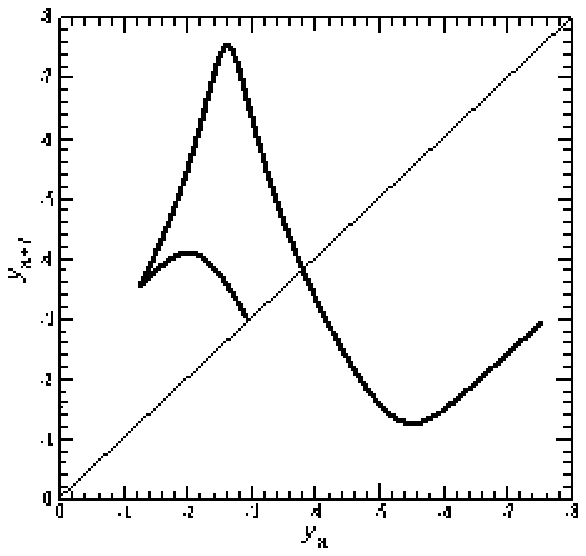} \\[-0.2cm]
      (a) Chaotic attractor & (b) First-return map \\[-0.2cm]
    \end{tabular}
    \caption{Sprott system R with a single fixed point. The trajectory wraps around component
$\phi_t$. One of the elliptic parabolo\"ids emerging from the fixed point results from a singularity
occuring when solving $\phi (x,y,z)=0$ as for Sprott system K. This leads to a funnel chaotic
attractor. Parameter values: $a=0.90$ and $b=0.395$.}
    \label{sprottR}
  \end{center}
\end{figure}

\subsection{Inverted Rössler-like chaos}

\begin{figure}[ht]
  \begin{center}
    \begin{tabular}{ccc}
      \includegraphics[height=4.0cm]{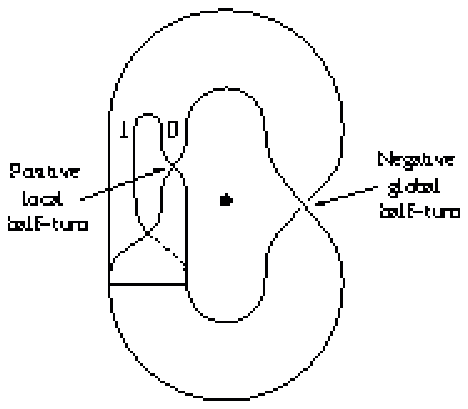} & &
      \includegraphics[height=4.0cm]{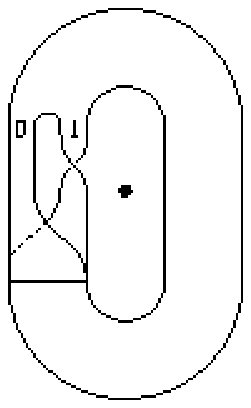} \\
      (a) Template with global torsion & & (b) Equivalent template \\[-0.2cm]
    \end{tabular}
    \caption{A template with a negative global half-turn and a positive
local half-turn (a) can be reduced under an isotopy to an inverted Rössler-like
template (b), that is, without any global half-turn and with a single local
half-turn (here negative). The branch with a local half-turn of the reduced
template is associated with the decreasing branch of the first-return map and
is located near the inner fixed point (designated by $\bullet$).}
    \label{gloeq}
  \end{center}
\end{figure}

Among Sprott systems with a single fixed point, two of them, namely systems L and N, produce
chaotic attractors which have an inverted Rössler-like topology. Typically, an inverted
Rössler-like attractor --- also named inverted Horseshoe attractor \cite{Gil02} --- differs
from a ``direct'' Rössler-like attractor by a global torsion of a half-turn. The usual organization
with the order preserving branch close to the inner fixed point and the order reversing branch
at the periphery of the attractor (Fig.\ \ref{templates}a) is therefore inverted and the
order reversing branch of the first-return map is close to the inner fixed point and the
order-preserving branch is at the periphery of the attractor (Fig.\ \ref{gloeq}b).

\begin{figure}[htbp]
    \begin{center}
      \begin{tabular}{cc}
        \includegraphics[height=5.5cm]{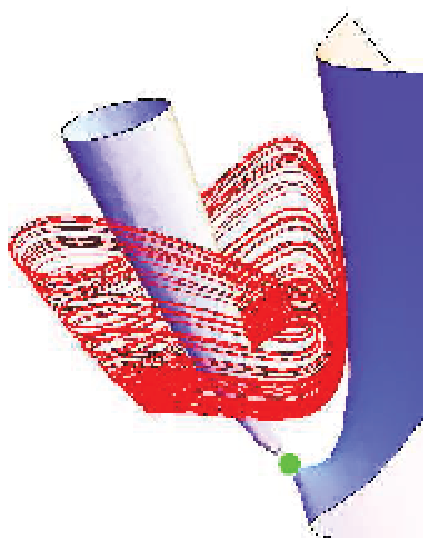} &
        \includegraphics[height=5.5cm]{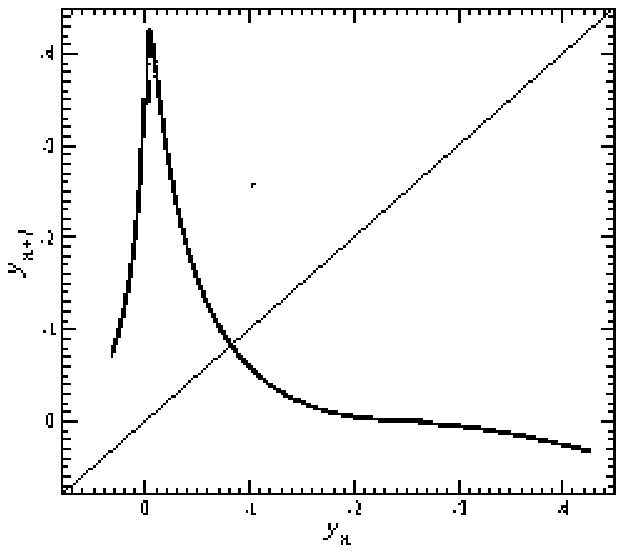} \\[-0.2cm]
        \multicolumn{2}{c}{(a) Sprott system L, $a=3.87$ and $b=0.91$} \\[0.2cm]
        \includegraphics[height=5.5cm]{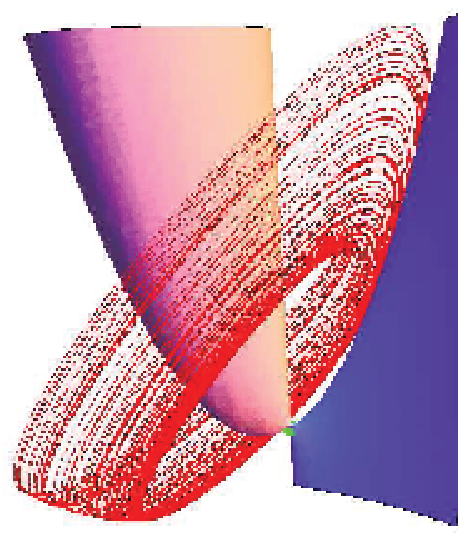}  &
        \includegraphics[height=5.5cm]{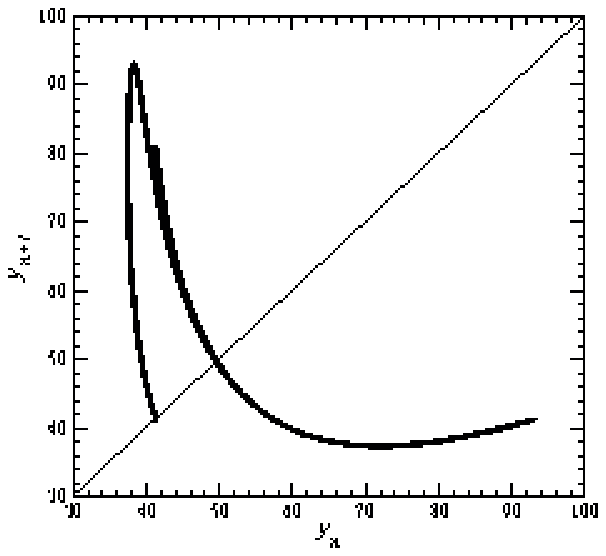}  \\[-0.2cm]
        \multicolumn{2}{c}{(b) Sprott system N, $a=4.2$} \\[-0.2cm]
      \end{tabular}
      \caption{Chaotic attractors solution to Sprott systems L and N. In both cases, the trajectory
spirals arouund component $\phi_t$ of the flow curvature manifold.}
      \label{sprottL}
    \end{center}
  \end{figure}

Sprott systems L and N do not present a component $\phi_t$ very different from those obtained for
systems I and J, for instance. Nevertheless, the two attractors (Figs.\ \ref{sprottL}a and
\ref{sprottL}b) are located relatively far from the fixed point (compared to previous cases).
In these two cases, the influence of component $\phi_t$ seems to
be induced by the second elliptic parabolo\"id which constrains the attractor by its periphery. Funnel
attractors would not be observed due to this external constraint.

\begin{figure}[ht]
  \begin{center}
    \begin{tabular}{cc}
      \includegraphics[height=5.3cm]{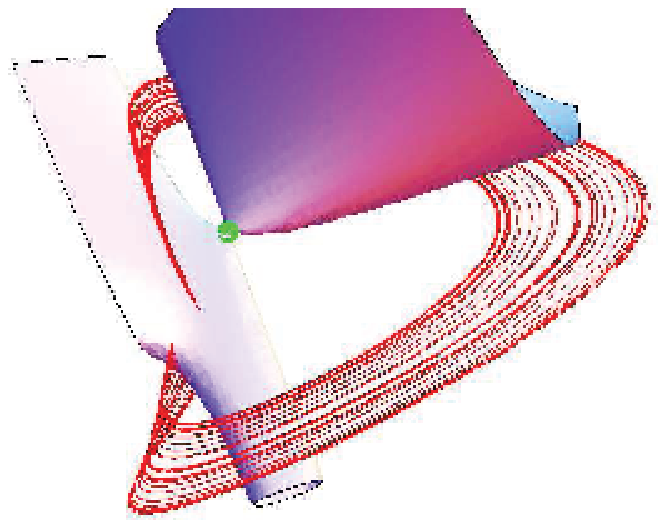} &
      \includegraphics[height=5.3cm]{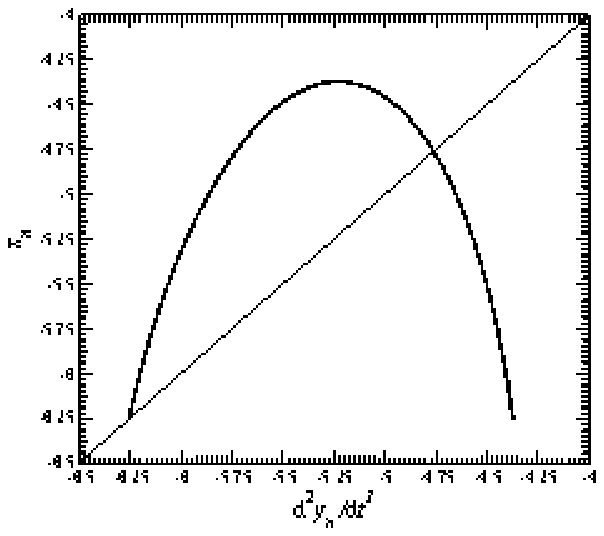} \\[-0.2cm]
      (a) Chaotic attractor & (b) First-return map \\[-0.2cm]
    \end{tabular}
    \caption{Chaotic attractor solution to Malasoma system A. Component $\phi_t$ presents an
unusual shape and crosses all the attractor, thus limiting the chaotic regime. Parameter value:
$a=2.017$. Initial conditions: $x_0=0.1$, $y_0=1$ and $z_0=1.9$.}
    \label{malA}
  \end{center}
\end{figure}

Two other systems with a single fixed point were proposed by Malasoma \cite{Mal02}. They are
minimal in the sense that it is not possible to obtain chaotic system with a simpler algebraic
structure. From the flow curvature manifold point of view, these two systems are similar and only
one of them is discussed here. As for many minimal systems, the chaotic domain in the parameter space
is quite limited. The attraction basin is also quite small. Component $\phi_t$ presents an unusual
shape with a cylindrical aspect for one of the two elliptic parabolo\"ids (Fig.\ \ref{malA}). Once again,
this results from numerical artefacts induced by a singularity appearing when solving $\phi
(x,y,z)=0$. In this case, the trajectory solution to Malasoma system A intersects component
$\phi_t$. Compared to all cases previously discussed, this is the first example for which the whole
attractor intersects component $\phi_t$ in the non ambiguous part. According to our assumption, such
a global intersection strongly limits the
development of the chaotic attractor. But the limitation of the dynamics occurs in a slighly
different way than the previous two cases. The first-return map presents a fully developed
unimodal map (Fig.\ \ref{malA}b), that is, more developed than those computed for Sprott systems
D and I (Figs.\ \ref{sprottD}). Nevertheless, real parts of the complex conjugated eigenvalues of
Malasoma system A is clearly non zero.  The intersection of the whole attractor with component
$\phi_t$ limits the region of the phase space where the attractor can exist. In particular, it
constrains the attractor to be developed quite far from the fixed point. As a consequence, the
branch without any half-turn is not observed and this is an inverted R\"ossler-like chaos.

Among these seven systems with a single fixed point, one --- Sprott system R --- presents a
time-dependent component $\phi_t$ around which the trajectory wraps (Fig.\ \ref{sprottR}).
Nevertheless, its time-dependent component is affected by a singularity which prevents us from avoiding
a spurious third elliptic parabolo\"id. It is therefore difficult to make conclusions about this system. The
presence
of a second fixed point is therefore not required to observe a chaotic attractor of a funnel
type. The relevant ingredient is indeed that the trajectory wraps around component $\phi_t$
without any intersection with it.

\section{Conclusion}
\label{conc}

It is still a very challenging problem to connect topological properties of phase portraits
with some analytical properties of the governing equations. Fixed points are certainly the first
step for such a connection. But the whole topological structure cannot be obtained from them.
In this paper, we showed that the flow curvature manifold can bring some additional light
on what structures the phase portrait. This manifold was split into one
time-dependent and one time-independent components. We showed that the time-independent component
was tangent to the osculating plane in the neighborhood of the inner fixed point. Our results
suggest that the time-dependent component is mainly responsible for limiting the development of
chaotic attractors when they are crossed by the trajectory. An attractor is thus not only
constrained by fixed points and some other solutions --- unstable periodic orbits for instance ---
co-existing in the phase space, but by the flow curvature manifold too. The next step is now
to investigate permeability properties of the flow curvature manifold to better understand why
the time dependent component $\phi_t$ of the flow curvature is not always crossed by
trajectories.

\paragraph*{Acknowledgements}
C. Letellier thanks L. A. Aguirre, R. Gilmore, U. Freitas and J.-M. Malasoma for stimulating
discussions.  Both of us thank Aziz-Alaoui for stimulating remarks while he was preparating his
own slides at the International Workshop-School {\it Chaos and Dynamics in Biological Networks} in
Carg\`ese (Corsica).

\end{document}